\documentclass[10pt,twosided]{article}
\usepackage[tbtags]{amsmath}
\allowdisplaybreaks[4]
\usepackage{epsfig}
\usepackage{graphicx, subfigure}
\usepackage{amssymb,color}

\definecolor{c20}{rgb}{0.,0.7,0.}
\definecolor{c30}{rgb}{0.,0.,1.}
\definecolor{c40}{rgb}{1,0.1,0.7}
\definecolor{c50}{rgb}{1,0,0}
\definecolor{c60}{rgb}{1,0.9,0.1}
\definecolor{c70}{rgb}{0.50,1.00,0.00}
\definecolor{c80}{rgb}{0.00,1.00,0.00}
\definecolor{c10}{rgb}{1.00,0.00,0.50}

\def\Ea#1{{\textcolor{c50}{#1}}}

\def\Ea#1{#1}

\def\EH#1{{\textcolor{c50}{#1}}}
\def\EH#1{#1}

\def\wz#1{{\textcolor{c30}{#1}}}
\def\wz#1{#1}
\def\zw#1{{\textcolor{c30}{#1}}}
\def\zw#1{#1}
\def\zc#1{{\textcolor{c30}{#1}}}
\def\zc#1{#1}

\def\Ha#1{{\textcolor{c60}{#1}}}
\def\Ha#1{#1}

\def\aH#1{{\textcolor{c60}{#1}}}
\def\aH#1{#1}

\def\eH#1{{\textcolor{c40}{#1}}}
\def\eH#1{#1}

\def\peng#1{{\textcolor{c80}{#1}}}
\def\peng#1{#1}

\def\pzx#1{{\textcolor{c80}{#1}}}
\def\pzx#1{#1}

\def\PZX#1{{\textcolor{c10}{#1}}}
\def\PZX#1{#1}
\newcommand{\pk}[1]{\mathbb{P} \left( #1 \right) }

\newcommand{\BQN}{\begin{eqnarray}}
\newcommand{\EQN}{\end{eqnarray}}
\newcommand{\BQNY}{\begin{eqnarray*}}
\newcommand{\EQNY}{\end{eqnarray*}}

\newcommand{\BS}{\begin{sat}}
\newcommand{\ES}{\end{sat}}
\newcommand{\BT}{\begin{theo}}
\newcommand{\ET}{\end{theo}}
\newcommand{\BK}{\begin{korr}}
\newcommand{\EK}{\end{korr}}

\newcommand{\BD}{\begin{de}}
\newcommand{\ED}{\end{de}}

\newcommand{\BIT}{\begin{itemize}}
\newcommand{\EIT}{\end{itemize}}
\newcommand{\BDI}{\begin{description}}
\newcommand{\EDI}{\end{description}}

\newcommand{\BRM}{\begin{remark}}
\newcommand{\ERM}{\end{remark}}

\newcommand{\BEL}{\begin{lem}}
\newcommand{\EEL}{\end{lem}}

\newtheorem{theo}{Theorem}[section]
\newtheorem{sat}[theo]{Proposition}
\newtheorem{de}[theo]{Definition}
\newtheorem{lem}{Lemma}[section]

\newtheorem{korr}[theo]{Corollary}
\newtheorem{remark}[theo]{Remark}

\newtheorem{theorem}{Theorem}[section]
\newtheorem{corollary}{Corollary}[section]
\newtheorem{lemma}{Lemma}[section]
\numberwithin{equation}{section}

\newcommand{\prooftheo}[1]{ \textsc{Proof of Theorem} \ref{#1} }

\newcommand{\prooflem}[1]{\textsc{Proof of Lemma} \ref{#1}}
\newcommand{\proofkorr}[1]{\textsc{Proof of Corollary} \ref{#1}}

\newcommand{\QED}{\hfill $\Box$}

\def\d{\mathrm{d}}

\newcommand{\nwc}{\newcommand}
\nwc{\COM}[1]{}

\def\IF{\infty}

\newcommand{\ABs}[1]{ \biggl \lvert #1 \biggr \rvert}

\def\DFxy{\Delta(F^{n}_{\rho_n}, H_\lambda;x,y)}
\def\DFxyA{\Delta(F^{n}_1, H_0;x,y)}
\def\DFxyB{\Delta(F^{n}_0, H_\IF;x,y)}

\def\DFxyD{\Delta(F^{n}_{\rho_n}, H_\IF;x,y)}
\def\DFxyE{\Delta(F^{n}_{\rho_n}, H_0;x,y)}
\begin{document}
\title{Higher-order expansions of distributions \\of maxima in a H\"{u}sler-Reiss model}
\author{{$^a$Enkelejd Hashorva \quad $^b$Zuoxiang Peng \quad $^a$Zhichao Weng}\\
{\small\it $^a$Faculty of Business and Economics (HEC), University of Lausanne, 1015 Lausanne, Switzerland} \\
{\small\it $^b$School of Mathematics and Statistics, Southwest
University, 400715 Chongqing, China } }
\maketitle

\begin{quote}
{\bf Abstract:} The max-stable H\"usler-Reiss distribution which
arises as the limit distribution of maxima of bivariate
Gaussian triangular arrays has been shown to be useful in various
extreme value models. For such triangular arrays, this paper
establishes higher-order asymptotic expansions of the joint distribution of
maxima under refined H\"{u}sler-Reiss conditions. In particular, the rate of convergence of normalized maxima to
\Ea{the} H\"usler-Reiss distribution is explicitly calculated. 
%Our findings are  supported by the results of a numerical analysis.

{\bf Key Words:} H\"{u}sler-Reiss max-stable distribution;
higher-order asymptotic expansion; triangular arrays; Gaussian random vector.
\end{quote}

\section{Introduction}
The fact that the componentwise maxima of bivariate \Ha{Gaussian
random vectors \Ea{possess} asymptotic} independent components (see
e.g., \cite{Faletal2010}) has been seen as a drawback in extreme
value theory since for modeling asymptotically dependent risks the
classical and tractable Gaussian framework is \peng{in}adequate. In
the seminal paper \cite{MR980699} this drawback was removed by
considering triangular arrays where the dependence \peng{may}
increase with $n$. \aH{Specifically}, \peng{let
$\{(X_{nk},Y_{nk}),1\le k\le n, n\ge 1\}$ be a triangular array of}
independent \aH{standard (mean-zero and unit variance)} bivariate
Gaussian random vectors \Ha{with correlations $\{\rho_n, n\ge 1\}$
and joint distribution \aH{function} \peng{$F_{\rho_n}$}}. The
principal finding of \cite{MR980699} is
\def\ux{u_n(x)}
\def\uy{u_n(y)}
\BQN\label{eqA} \lim_{n\to \IF} \sup_{x,y\in \mathbb{R}}\ABs{
F_{\rho_n}^n( x/b_n+ b_n,\; y/b_n+ b_n)-H_{\lambda}(x,y)}=0,
\EQN
provided that the so-called H\"usler-Reiss condition
\begin{equation}\label{eq1.1}
\lim_{n\to\IF}\frac{1}{2}b_n^2(1-\rho_n)=\lambda^2\peng{\quad\mbox{with}\;
\lambda}\in[0,\IF]
\end{equation}
holds \Ha{with} $b_n$ given by
\begin{equation}\label{eqab}
n(1-\Phi(b_n))=1
\end{equation}
or $nb_{n}^{-1}\varphi(b_{n})=1$, where \Ha{$\Phi$ denotes} the $N(0,1)$ distribution function and $\varphi(x)=\Phi^{\prime}(x)$;
\aH{see
 \cite{MR980699} for more details.} The max-stable H\"usler-Reiss
distribution $H_\lambda$ is given by
\[H_{\lambda}(x,y)=\exp\left(-\Phi\left(\lambda+\frac{y-x}{2\lambda}\right)e^{-x}-\Phi\left(\lambda+\frac{x-y}{2\lambda}\right)e^{-y}\right),
\quad \aH{x,y\in \mathbb{R}},
\]
\Ha{with  $H_0(x,y)=\exp(-e^{-\min(x,y)})$} and
$H_{\IF}(x,y)=\Lambda(x)\Lambda(y)$, \eH{where
\peng{$\Lambda(x)=\exp(-e^{-x}),x\in \mathbb{R}$} \Ha{is} the
Gumbel distribution.}

 \eH{In fact, the bivariate H\"usler-Reiss
distribution appeared in another context in \cite{bro1977}, see for
recent contribution in this direction \peng{\cite{eng2011, kab2011, oes2012,Das14}}. Related results for more general triangular
arrays can be found in \peng{\cite{H08Dirichlet, eng2012b, eng2012a, MR2137118, HashKab, MR2998759, MR1421762, MR1325215, HB, Deb14, Weng14}};
 \Ea{an interesting} statistical applications related to \Ea{the} H\"usler-Reiss distribution \EH{is presented in} \cite{eng2012a}.\\
 For both applications and \EH{various} theoretical investigations\peng{,} it is of interest to know how good \Ea{the} H\"usler-Reiss distribution approximates the distribution of the
bivariate maxima. \peng{So,} \aH{a} natural goal of this paper is to
investigate the rate of convergence in \eqref{eqA}, \peng{i.e., the}
speed of convergence to $0$ \EH{as $n\to \IF$} of the following difference
$$ \DFxy:= F_{\rho_n}^n(\ux,\uy)-H_{\lambda}(x,y),$$
\PZX{where $\Ea{u_n}(s)= s/b_n+ b_n$ with norming constant $b_{n}$ given by
\eqref{eqab}.} \aH{In the literature the only available results
\Ha{concern} the univariate problem, namely in \cite{MR606804} it
has been shown that}
\begin{equation}\label{eq1.3}
\lim_{n\to
\IF}b_n^2\Big[b_n^2\Big(\Phi^n(\ux)-\Lambda(x)\Big)-s(x)\Lambda(x)\Big]=\left(t(x)+\frac{1}{2}s^2(x)\right)\Lambda(x),
\end{equation}
with  $b_n$ given by \eqref{eqab} and $s(x),t(x)$ defined as
\begin{equation}\label{eqkl}
s(x)=2^{-1}(x^2+2x)e^{-x} \quad \mbox{  and  }\quad t(x)=-8^{-1}(x^4+4x^3+8x^2+16x)e^{-x}.
\end{equation}
}

\eH{\peng{In order to derive the rate of convergence of \aH{$\DFxy$ to 0}
\aH{we shall} introduce a refinement of the H\"usler-Reiss condition \eqref{eq1.1}, namely we
shall suppose that
\begin{eqnarray}\label{eq2.1}
\lim_{n \to \IF} b^2_n(\lambda-\lambda_n)=\alpha\in \mathbb{R}
\end{eqnarray}
holds with \peng{$\lambda_n=(b_n^2(1-\rho_n)/2)^{1/2}$ and
$\lambda\in (0,\infty)$}.  \aH{By assuming further that}
$\delta_n=b^2_n(\lambda-\lambda_n)-\alpha$ also converges to $0$
with a speed \aH{determined} again by $b_n^2$, we are able to refine \peng{the second-order
approximation significantly. The analysis of the two extreme cases
$\lambda=0$ and $\lambda=\IF$ are more complicated and more
information related to $\rho_{n}$ is needed. \Ha{Two special cases}  $\rho_n=1$ and $\rho_n \in [-1, 0]$ for all large $n$ are
\aH{explicitly} solved}. %A numerical study \aH{presented in the next section} shows the practical value of our \peng{higher-order} approximations.
}}

\peng{The rest of the paper is organized as follows.} In Section
\ref{sec2} we present the main results .
%\EH{and then illustrate our} findings \peng{by a small numerical  \Ea{study}}. 
All the proofs are relegated to Section \ref{sec3}.

\section{Main Results}\label{sec2}
\eH{In the following we shall denote throughout by $b_n$ the constants defined in \eqref{eqab} and
%write $u_n(z)=b_n+z/b_n, z\in \mathbb{R}$.
further $\lambda$ shall always be defined with respect \aH{to the H\"usler-Reiss condition}
\eqref{eq1.1}.} \peng{\aH{Next}, we \Ea{derive}} the \peng{second-order
expansions of bivariate extremes} under the second-order
H\"usler-Reiss condition \eqref{eq2.1}.

\def\kapx{\kappa(\alpha, \lambda, x,y)}
\def\tapx{\tau(\alpha, \beta, \lambda, x,y)}

\begin{theorem}\label{th2.1}
\eH{If \eqref{eq2.1}
holds} with \peng{$\lambda_n=(b_n^2(1-\rho_n)/2)^{1/2}$} and
$\lambda\in(0,\IF)$, then for all $x,y \in \mathbb{R}$ we have
\BQN\label{r:thA}
\lim_{n\to \IF}b_n^2 \DFxy %\Bigl( F_{\rho_n}^n(u_n(x),u_n(y))-H_{\lambda}(x,y)\Bigr)
&=&\kapx H_{\lambda}(x,y), \EQN
where
\begin{eqnarray*}
\kapx&=&\eH{s(x)}\Phi\left(\lambda+\frac{y-x}{2\lambda}\right)
+\eH{s(y)}\Phi\left(\lambda+\frac{x-y}{2\lambda}\right)%\\
%&&
+(2\alpha-\peng{\lambda(\lambda^{2}+x+
y+2)})e^{-x}\varphi\left(\lambda+\frac{y-x}{2\lambda}\right),
\end{eqnarray*}
\peng{where $s(z), z\in \mathbb{R}$ is defined by \eqref{eqkl}.}
\end{theorem}

If the second-order H\"usler-Reiss condition is further refined
\peng{to a third-order one,} a finer result \peng{than} that stated
in \eqref{r:thA} can be obtained. Indeed, this can be achieved by
\Ea{introducing} a restriction on the difference $\delta_n:=
b^2_n(\lambda-\lambda_n)-\alpha$, \peng{namely}
\begin{eqnarray}\label{eq2.2}
\lim_{n \to \IF}b^2_n \delta_n =\beta\in \mathbb{R}.
\end{eqnarray}
\peng{\aH{Utilising further} condition \eqref{eq2.2} we \aH{derive below} \EH{a}
{third-order expansion} of \EH{the} joint \aH{distribution} of extremes. \Ea{For simplicity we shall omit the expression of the function $\tau$ below, it is specified in \eqref{addpeng3}}.}
\begin{theorem}\label{th2.2}
If \eqref{eq2.2} holds \peng{with $\lambda\in(0,\IF)$}, then for all
$x,y \in \mathbb{R}$ we have
\begin{equation}\label{r:th2}
\lim_{n\to \IF}b_n^2 \Biggl(b_n^2
\DFxy- \kapx H_{\lambda}(x,y)\Biggr)
=\left( \tapx +\frac{1}{2}\kappa^2(\alpha,\lambda, x,y)\right)H_{\lambda}(x,y).
\end{equation}
%\peng{where $\tapx $ is given by
%\eqref{addpeng3}.}
\end{theorem}

For the two extreme cases $\lambda=0$ and $\lambda=\infty$ we
first consider \pzx{two special cases satisfied for all large $n$, \EH{namely}  $\rho_{n}\in [-1,0]$ and
$\rho_{n}= 1$  including components of each Gaussian vector
with independence  ($\rho_{n}=0$), complete negative
dependence ($\rho_{n}= -1$)  and complete positive dependence
($\rho_{n}=1$), respectively.} \pzx{
\begin{theorem}\label{th2.3}
Let $s(z)$ and $t(z)$ be those defined as in \eqref{eqkl} and set
$u_n(z)=b_n+z/b_n$, $z\in \mathbb{R}$.
\begin{itemize}
\item[(i).] For $\rho_{n}\in [-1,0],n\ge 1$ and any $x,y \in
\mathbb{R}$ we have
\begin{eqnarray}\label{eq2.5}
\lim_{n\to \IF} b_n^2\Big[b_n^2 \DFxyD-(s(x)+s(y))H_{\IF}(x,y)\Big]
=\left(t(x)+t(y)+\frac{1}{2}(s(x)+s(y))^2\right)H_{\IF}(x,y).
\end{eqnarray}
\item[(ii).] \EH{If}  $\rho_{n}=1, n\ge 1$, \EH{then} for any $x,y \in
\mathbb{R}$ we have
\begin{eqnarray}\label{eq2.4}
\lim_{n\to \IF}b_n^2\Big[b_n^2 \DFxyA -s(\min(x,y))H_{0}(x,y)\Big]
=\left(t(\min(x,y))+\frac{1}{2}(s(\min(x,y)))^2\right)H_{0}(x,y).
\end{eqnarray}
\end{itemize}
\end{theorem}}

\aH{We consider next the other cases} of $\rho_{n}\in (0,1)$ such that
$\lambda_{n}\to 0$ or $\lambda_{n}\to \infty$. With more information
on the asymptotic behavior of $\rho_{n}$ \aH{we obtain below upper bounds for the}
convergence rates \Ha{of} $F_{\rho_{n}}^{n}$ to \aH{$H_0$ or $H_\IF$}.

\begin{corollary}\label{cor2.1}
For some $\mathbb{C}>0$ and $R(x,y)=\mathbb{C}(\exp(2|x|)+\exp(2|y|))$ we have:\\ %Let the norming constant $b_{n}$ be given by \eqref{eqab} and $s(z)$
%is the one defined in \eqref{eqkl}. For $\rho_{n}\in(0,1)$,
(i). Suppose that  $\rho_n\in (0,1), n\ge 1$ and \eqref{eq1.1} holds with $\lambda=\IF$. \aH{If further}
\peng{$\zw{\frac{1}{2}}((1-\rho_{n})\ln  n -(2+\rho_{n})\ln \ln  n)\to\EH{\gamma \in (-\IF,\IF]}$ as
$n\to\infty$}, then for all $x,y \in \mathbb{R}$
\begin{eqnarray*}
\limsup_{n\to \IF}b_n^2\ABs{\DFxyD }\le (|s(x)|+|s(y)|)H_{\IF}(x,y) + \EH{ e^{-\gamma} R(x,y)}.
\end{eqnarray*}
(ii). \aH{If}  \wz{$(1-\rho_n)(\ln  n)^3 \to \EH{\tau^2} \in [0,\IF)$}
as $n\to\infty$, then for all $x,y \in \mathbb{R}$
\begin{eqnarray*}
\limsup_{n\to \IF}b_n^2\ABs{\DFxyE }\le |s(\min(x,y))|H_{0}(x,y)+
\PZX{\tau}R(x,y).
\end{eqnarray*}
\end{corollary}

\peng{
\begin{remark}
For \EH{the} H\"{u}sler-Reiss model the rates of convergence of
$F_{\rho_{n}}^{n}(u_n(x),u_n(y))$ to its ultimate
max-stable distribution $H_{\lambda}(x,y)$ is  proportional to
$O(1/\ln  n)$ for all cases studied in this paper.
\end{remark}}

\COM{
\peng{To this end, \aH{we shall present} a small numerical study which
illustrates our findings, i.e., the accuracy of higher-order
expansions of $F_{\rho_{n}}^{n}$ by \aH{the corresponding} max-stable H\"{u}sler-Reiss
distribution. \aH{We shall \zw{discuss} three particular cases:} \\
 \aH{ $(i)$} $\rho_{n}\in
[-1,0]$ implying $\lambda=\infty$;\\
\aH{ $(ii)$} $\rho_{n}=1$ implying
$\lambda=0$;\\
\aH{$(iii)$} $\lambda\in (0,\infty)$ with}
\begin{equation}\label{eq4.1}
\rho_n=1-\frac{2\lambda^2}{b_n^2}+\frac{4\alpha
\lambda}{b_n^4}+\frac{4\beta\lambda-2\alpha^{2}}{b_n^6},
\end{equation}
\peng{where $b_n$ \aH{satisfies} \eqref{eqab}, which implies that  condition \eqref{eq2.2} holds. For finite $n$, we calculate the
actual values $F^n_{\rho_n}(u_n(x),u_n(y))$, the first-order
asymptotics $H_{\lambda}(x,y)$, the second-order asymptotics and the
third-order asymptotics according to the values of $\rho_{n}$, i.e.,
\begin{itemize}
\item[ $(i)$.] if $\rho_{n}\in [-1,0]$, by \eqref{eq2.5} the second-order
and the third-order asymptotics are respectively given by
$H_{\IF}(x,y)\Big[1+\frac{s(x)+s(y)}{b_{n}^{2}}\Big]$ and
$H_{\IF}(x,y)\Big[1+\frac{s(x)+s(y)}{b_{n}^{2}}+\frac{t(x)+t(y)+(s(x)+s(y))^{2}/2}{b_{n}^{4}}\Big]$;
\item[$(ii)$.] if $\rho_{n}=1,n\ge 1$ by \eqref{eq2.4}  the second-order
and the third-order asymptotics are given by
$H_{0}(x,y)\Big[1+\frac{s(\min(x,y))}{b_{n}^{2}}\Big]$ and
$H_{0}(x,y)\Big[1+\frac{s(\min(x,y))}{b_{n}^{2}}+\frac{t(\min(x,y))+(s(\min(x,y)))^{2}/2}{b_{n}^{4}}\Big]$,
respectively; and
\item[$(iii)$.] if $\rho_{n}$ is given by \eqref{eq4.1} with
fixed $\lambda$, $\alpha$ and $\beta$, \EH{then in view of} \eqref{r:th2}  the
second-order and the third-order asymptotics are given by
$H_{\lambda}(x,y)(1+b_n^{-2} \kapx )$ and
$H_{\lambda}(x,y)[1+b_n^{-2}\kapx +b_n^{-4}(\tapx +2^{-1}(\kapx )^2)]$,
respectively.
\end{itemize}}
}
\COM{ 
Figure 1 compares the actual values with above three asymptotics
\wz{with $x=y$}. Figure 2 compares the difference of the actual
value with above three asymptotics by contour line in the plane.
According to \peng{Figures} 1 and 2, we observe that the third-order
asymptotics is closest to the actual
values.

\begin{figure}
\begin{center}
\subfigure[$\rho_n= -1$]
{%
 \epsfig{file=rho-1.eps, height=120pt, width=200pt,angle=0}
                   }%
\subfigure[$\lambda=2.5,\alpha=-5,\beta=10, \rho_n=-0.80$]
{%
 \epsfig{file=rho-0.80.eps, height=120pt, width=200pt,angle=0}
                   }%
                   \\
\subfigure[$\lambda=2.5,\alpha=-2,\beta=5, \rho_n=-0.48$]
{%
 \epsfig{file=rho-0.48.eps, height=120pt, width=200pt,angle=0}
                   }%
 \subfigure[$\rho_n= 0$]
{%
 \epsfig{file=rho0.eps, height=120pt, width=200pt,angle=0}
                   }%
                   \\
 \subfigure[$\lambda=2,\alpha=2,\beta=-10, \rho_n=0.23$]
{%
 \epsfig{file=rho0.23.eps, height=120pt, width=200pt,angle=0}
                   }%
\subfigure[$ \lambda=2,\alpha=3,\beta=11,\rho_n=0.51$]
{%
\epsfig{file=rho0.51.eps, height=120pt, width=200pt,angle=0}
                   }%
                   \\
 \subfigure[$ \lambda=1,\alpha=2,\beta=5,\rho_n=0.89$]
{%
\epsfig{file=rho0.89.eps, height=120pt, width=200pt,angle=0}
                   }%
\subfigure[$ \rho_n=1$]
{%
\epsfig{file=rho1.eps, height=120pt, width=200pt,angle=0}
                   }%
\caption{\peng{Actual values and its approximations with $n=10^{3},
x=y\in [-2,8]$. The actual values with black color, the first-order
asymptotics with blue color, the second-order asymptotics with green
color and the third-order asymptotics with red color.}}
\end{center}
\end{figure}

\begin{figure}
\begin{center}
\subfigure[$\rho_n= -1$]
{%
 \epsfig{file=crho-1.eps, height=120pt, width=200pt,angle=0}
                   }%
\subfigure[$\lambda=2.5,\alpha=-5,\beta=10, \rho_n=-0.80$]
{%
 \epsfig{file=crho-0.80.eps, height=120pt, width=200pt,angle=0}
                   }%
                   \\
\subfigure[$\lambda=2.5,\alpha=-2,\beta=5, \rho_n=-0.48$]
{%
 \epsfig{file=crho-0.48.eps, height=120pt, width=200pt,angle=0}
                   }%
 \subfigure[$\rho_n= 0$]
{%
 \epsfig{file=crho0.eps, height=120pt, width=200pt,angle=0}
                   }%
                   \\
 \subfigure[$\lambda=2,\alpha=2,\beta=-10, \rho_n=0.23$]
{%
 \epsfig{file=crho0.23.eps, height=120pt, width=200pt,angle=0}
                   }%
\subfigure[$ \lambda=2,\alpha=3,\beta=11,\rho_n=0.51$]
{%
\epsfig{file=crho0.51.eps, height=120pt, width=200pt,angle=0}
                   }%
                   \\
 \subfigure[$ \lambda=1,\alpha=2,\beta=5,\rho_n=0.89$]
{%
\epsfig{file=crho0.89.eps, height=120pt, width=200pt,angle=0}
                   }%
\subfigure[$ \rho_n=1$]
{%
\epsfig{file=crho1.eps, height=120pt, width=200pt,angle=0}
                   }%
\caption{\wz{The contour line of} actual values and its approximations with $n=10^{3},
x=y\in [-2,8]$. The actual values with black color, the first-order
asymptotics with blue color, the second-order asymptotics with green
color and the third-order asymptotics with red color.}
\end{center}
\end{figure}
}

\def\pHn{\overline{\Phi}_n}

\section{Proofs}\label{sec3}
\eH{Recall that we set
$u_n(x)=b_n+x/b_n, x\in \mathbb{R}$ \wz{with $b_n$
satisfying equation \eqref{eqab}}. Define further below
$$ \overline\Phi(x)= 1- \Phi(x), \quad \pHn(s)= n \overline\Phi(u_n(s))$$
and
$$\EH{I_k:=}\int_y^{\IF}\varphi\left(\lambda+\frac{x-z}{2\lambda}\right)e^{-z}z^k \d z, \quad k=0,\cdots, 3.$$}
The following formulas obtained by partial integration will be used in the proofs below:
\BQN\label{eq3.10}
I_0&=&2\lambda e^{-x}\overline \Phi\left(\lambda+\frac{y-x}{2\lambda}\right),\\
\label{eq3.11}
I_1&=&(2\lambda x-4\lambda^3) e^{-x}\overline \Phi\left(\lambda+\frac{y-x}{2\lambda}\right)
+4\lambda^2 e^{-x}\varphi\left(\lambda+\frac{y-x}{2\lambda}\right),\\
I_2&=&(8\lambda^5-8\lambda^3 x+8\lambda^3+2\lambda x^2)e^{-x}\overline \Phi\left(\lambda+\frac{y-x}{2\lambda}\right)
+(-8\lambda^4+4\lambda^2 x+4\lambda^2 y)e^{-x}\varphi\left(\lambda+\frac{y-x}{2\lambda}\right)
\label{I2}\\
I_3&=&(24\lambda^5 x-12\lambda^3x^2+24\lambda^3x+2\lambda x^3-16\lambda^7-48\lambda^5)e^{-x}
\overline \Phi\left(\lambda+\frac{y-x}{2\lambda}\right)
\notag\\
&&+(16\lambda^6-16\lambda^4x-8\lambda^4y+32\lambda^4+4\lambda^2x^2+4\lambda^2xy+4\lambda^2y^2)e^{-x}\varphi\left(\lambda+\frac{y-x}{2\lambda}\right).
\label{I3}
\EQN

\begin{lemma}\label{le3.1}
\aH{If} $(X,Y)$ is  a bivariate normal vector with correlation $\rho \eH{\in (-1,1)}$, then
\begin{eqnarray}
&&n\pk{X>u_n(x),Y>u_n(y)}\nonumber\\
&=&\pHn(y)-\int_y^{\IF}\Phi\left(\frac{u_n(x)-\rho
u_n(z)}{\sqrt{1-\rho^2}}\right)
e^{-z}\left[1+\left(1-\frac{z^2}{2}\right)\frac{1}{b_n^{2}}+\left(\frac{z^4}{8}-\frac{z^2}{2}-2\right)\frac{1}{b_n^{4}}\right] \d z +O(b_n^{-6})
\quad \quad \label{eq3.4} \\
&=&\pHn(y)-\int_y^{\IF}\Phi\left(\frac{u_n(x)-\rho
u_n(z)}{\sqrt{1-\rho^2}}\right)
e^{-z}\left[1+\left(1-\frac{z^2}{2}\right)\frac{1}{b_n^{2}}\right]
\d z +O(b_n^{-4}).\label{eq3.5}
\end{eqnarray}
\end{lemma}

\prooflem{le3.1} First note \aH{that}
$$\left|e^{-x}-\left(1-x+\frac{x^2}{2}\right)\right|<\zc{\frac{x^3}{6}+\frac{x^4}{24}}$$
for $x>0$, which \peng{implies}
\begin{eqnarray*}%\label{eq3.3}
\int_{u_n(y)}^{\IF}\Phi\left(\frac{u_n(x)-\rho z}{\sqrt{1-\rho^2}}\right)\varphi(z) \d z
&=&b_n^{-1}\varphi(b_n)\left[\int_y^{\IF}\Phi\left(\frac{u_n(x)-\rho
u_n(z)}{\sqrt{1-\rho^2}}\right)
e^{-z}\left(1-\frac{z^2}{2b_n^{2}}+\frac{z^4}{8b_n^{4}}\right) \d z
+O(b_n^{-6})\right]
\end{eqnarray*}
for large $n$. Hence
\begin{eqnarray*}
\lefteqn{\pk{X>u_n(x),Y>u_n(y)}}\notag\\
&=&\int_{u_n(y)}^{\IF} \overline \Phi\left(\frac{u_n(x)-\rho
z}{\sqrt{1-\rho^2}}\right)\varphi(z) \d z\nonumber\\
&=&\overline \Phi(u_n(y))-b_n^{-1}\varphi(b_n)\int_y^{\IF}\Phi\left(\frac{u_n(x)-\rho
u_n(z)}{\sqrt{1-\rho^2}}\right)
e^{-z}\left(1-\frac{z^2}{2b_n^{2}}+\frac{z^4}{8b_n^{4}}\right) \d z +O(b_n^{-7}\varphi(b_n)) \label{eq3.1}\\
&=&\overline \Phi(u_n(y))-b_n^{-1}\varphi(b_n)\int_y^{\IF}\Phi\left(\frac{u_n(x)-\rho
u_n(z)}{\sqrt{1-\rho^2}}\right)
e^{-z}\left(1-\frac{z^2}{2b_n^{2}}\right) \d z
+O(b_n^{-5}\varphi(b_n)) \label{eq3.2}.
\end{eqnarray*}
According to the definition of $b_n$ we have
$$n^{-1}=\overline \Phi(b_n)=b_n^{-1}\varphi(b_n)(1-b_n^{-2}+3b_n^{-4}+O(b_n^{-6}))$$
for large $n$ %, cf. \cite{MR912173},
thus the claim follows. \QED

For notational simplicity hereafter we set
\begin{eqnarray*}
A_{1n}=b_n^2\left(\lambda-\lambda_n\left(1-\frac{\lambda_n^2}{b_n^2}\right)^{-\frac{1}{2}}\right), \qquad
A_{2n}=\frac{1}{2}b_n^2\left(\frac{1}{\lambda}-\frac{1}{\lambda_n}\left(1-\frac{\lambda_n^2}{b_n^2}\right)^{-\frac{1}{2}}\right)
\end{eqnarray*}
and
$$A_{3n}=\lambda_n\left(1-\frac{\lambda_n^2}{b_n^2}\right)^{-\frac{1}{2}}.$$

\begin{lemma}\label{le3.2}
Under the conditions of Theorem \ref{th2.1}, we have
$$\lim_{n\to \IF}b_n^2\int_y^{\IF}\left(\Phi\left(\lambda+\frac{x-z}{2\lambda}\right)
-\Phi\left(\frac{u_n(x)-\rho_n
u_n(z)}{\sqrt{1-\rho_n^2}}\right)\right)e^{-z} \d z =
\kappa_1(\alpha,\lambda, x,y,),$$ where
$$\kappa_1(\alpha,\lambda,x,y)=(2\lambda^4-2\lambda^2x)e^{-x}\overline \Phi\left(\lambda+\frac{y-x}{2\lambda}\right)
+(2\alpha-3\lambda^3)e^{-x}\varphi\left(\lambda+\frac{y-x}{2\lambda}\right).$$
\end{lemma}

\prooflem{le3.2} Using the assumption \eqref{eq2.1} we have
\BQNY
\lim_{n\to \IF}A_{1n}&=&\lim_{n\to \IF}b_n^{2}\left(\lambda-\lambda_n-\frac{1}{2b_n^2}\lambda_n^3+O(b_n^{-4})\right)=\alpha-\frac{1}{2}\lambda^3,\label{eq3.6}\\
\lim_{n\to \IF}A_{2n}&=&\lim_{n\to \IF}\frac{1}{2}
b_n^{2}\left(\frac{\lambda_n-\lambda}{\lambda\lambda_n}-\frac{1}{2b_n^2}\lambda_n+O(b_n^{-4})\right)
=-\frac{1}{2}\alpha\lambda^{-2}-\frac{1}{4}\lambda,\\\label{eq3.7}
\label{eq3.8}
\lim_{n\to\IF}A_{3n}&=&\lim_{n\to\IF}\lambda_n(1+O(b_n^{-2}))=\lambda.
\EQNY
Hence since
\begin{equation}\label{eq3.9}
\frac{u_n(x)-\rho_n u_n(z)}{\sqrt{1-\rho_n^2}}
=\left(\lambda_n+\frac{x-z}{2\lambda_n}+\frac{\lambda_nz}{b_n^2}\right)
\left(1-\frac{\lambda_n^2}{b_n^2}\right)^{-\frac{1}{2}}\zc{\to \lambda+\frac{x-z}{2\lambda}}, \quad n\to \IF,
\end{equation}
\aH{then we obtain}
\begin{eqnarray}\label{eq3.12}
&&b_n^2\int_y^{\IF}\left(\lambda+\frac{x-z}{2\lambda}-\frac{u_n(x)-\rho_n u_n(z)}{\sqrt{1-\rho_n^2}}\right)
\varphi\left(\lambda+\frac{x-z}{2\lambda}\right)e^{-z} \d z \nonumber\\
&=&(A_{1n}+A_{2n}x) \EH{I_0} %\int_y^{\IF}\varphi\left(\lambda+\frac{x-z}{2\lambda}\right)e^{-z} \d z
-(A_{2n}+A_{3n})\EH{I_1}%\int_y^{\IF}\varphi\left(\lambda+\frac{x-z}{2\lambda}\right)e^{-z}z \d z
\nonumber\\
&\to&\left(\alpha-\frac{1}{2}\lambda^3-\frac{1}{2}\alpha  \lambda^{-2}x-\frac{1}{4} \lambda x\right) \EH{I_0}%\int_y^{\IF}\varphi\left(\lambda+\frac{x-z}{2\lambda}\right)e^{-z} \d z
-\left(\frac{3}{4}\lambda-\frac{1}{2}\alpha \lambda^{-2}\right)\EH{I_1}
%\int_y^{\IF}\varphi\left(\lambda+\frac{x-z}{2\lambda}\right)e^{-z}z \d z
\nonumber\\
&=&(2\lambda^4-2\lambda^2 x)e^{-x}\overline \Phi\left(\lambda+\frac{y-x}{2\lambda}\right)
+(2\alpha-3\lambda^3)e^{-x}\varphi\left(\lambda+\frac{y-x}{2\lambda}\right)
\end{eqnarray}
as $n\to \IF$. Using Taylor's expansion with Lagrange remainder term, we have
\peng{
\begin{eqnarray}\label{eq3.13}
\Phi\left(\frac{u_n(x)-\rho_n u_n(z)}{\sqrt{1-\rho_n^2}}\right)
&=&\Phi\left(\lambda+\frac{x-z}{2\lambda}\right)+
\varphi\left(\lambda+\frac{x-z}{2\lambda}\right)
v_n(x,y,\lambda) \notag \\%\left(\frac{u_n(x)-\rho_n u_n(z)}{\sqrt{1-\rho_n^2}}-\lambda-\frac{x-z}{2\lambda}\right)\nonumber\\
&&\qquad\qquad-\frac{1}{2}\EH{\xi_{n}(x,z)}\varphi(\EH{\xi_{n}(x,z)}) v_n^2(x,z,\lambda)
\end{eqnarray}
with $v_n(x,z,\lambda): =\left(\frac{u_n(x)-\rho_n
u_n(z)}{\sqrt{1-\rho_n^2}}-\lambda-\frac{x-z}{2\lambda}\right)$ and some $\EH{\xi_{n}(x,z)}$ between $\frac{u_n(x)-\rho_n
u_n(z)}{\sqrt{1-\rho_n^2}}$ and $\lambda+\frac{x-z}{2\lambda}$.
Moreover, \PZX{by arguments similar to \eqref{eq3.12},}
\zc{combining with \eqref{eq3.9}} we have
\begin{eqnarray*}
&&b_n^2\int_y^{\IF}v_n^2(x,z,\lambda)
\EH{\xi_{n}(x,z)}\varphi(\EH{\xi_{n}(x,z)})e^{-z}\d z\nonumber\\
&=&\zc{b_n^{-2}\int_y^{\IF}\left(A_{1n}+A_{2n}x-(A_{2n}+A_{3n})z\right)^2
\EH{\xi_{n}(x,z)}\varphi(\EH{\xi_{n}(x,z)})e^{-z}\d z}\nonumber\\
&=&O(b_n^{-2}),
\end{eqnarray*}}
which together with \eqref{eq3.12} and \eqref{eq3.13} established the proof. \QED

\begin{lemma}\label{le3.3}
Under the conditions of Theorem \ref{th2.2}, we have
$$\lim_{n\to \IF}b_n^2\left[b_n^2\int_y^{\IF}\left(\Phi\left(\lambda+\frac{x-z}{2\lambda}\right)
-\Phi\left(\frac{u_n(x)-\rho_n u_n(z)}{\sqrt{1-\rho_n^2}}\right)\right)e^{-z} \d z -\kappa_1(\alpha,\lambda,x,y)\right]
=\tau_1(\alpha,\beta,\lambda,x,y),$$
where
$\kappa_1(\alpha,\lambda,x,y)$ is defined in Lemma \ref{le3.2} and
\begin{eqnarray*}
\tau_1(\alpha,\beta,\lambda,x,y)&=&(2\lambda^8+8\lambda^6-4\lambda^6 x+2\lambda^4x^2-4\lambda^4x-8\alpha\lambda^3+4\alpha \lambda x)e^{-x}\overline \Phi\left(\lambda+\frac{y-x}{2\lambda}\right)\\
&+&\left(2\beta+9\alpha\lambda^2-\frac{23}{4}\lambda^5-\frac{3}{8}\lambda^3xy-\alpha
\lambda^2x+\frac{3}{4}\alpha y^2-\frac{1}{4}\alpha^2 \lambda^{-3}y^2
-\frac{1}{4}\alpha^2 \lambda^{-3}x^2-\alpha
\lambda^2y-\frac{1}{4}\alpha x^2
-\frac{7}{4}\lambda^7\right.\\
&+&\left.\frac{7}{2}\lambda^5x-\frac{1}{16}\lambda^3x^2
-\alpha\lambda^4+\alpha^2\lambda+\frac{3}{2}\lambda^5y-\frac{9}{16}\lambda^3y^2-\frac{1}{2}\alpha
xy+\frac{1}{2}\alpha^2\lambda^{-3}xy\right)
e^{-x}\varphi\left(\lambda+\frac{y-x}{2\lambda}\right).
\end{eqnarray*}
\end{lemma}

\prooflem{le3.3} \peng{By} assumption \eqref{eq2.2} we have
\begin{eqnarray*}\label{eq3.15}
\lim_{n\to \IF}b_n^2\left(A_{1n}-\alpha+\frac{1}{2}\lambda^3\right)
%&=&\lim_{n\to %\IF}b_n^2\left[b_n^2\left(\lambda-\lambda_n-\frac{\lambda_n^3}{2b_n^2}-\frac{3\lambda^5}{8b_n^4}+O(b_n^{-6})\right)-\alpha+\frac{1}{2}\lambda^3\right]\nonumber\\
&=&\beta+\frac{3}{2}\alpha\lambda^2-\frac{3}{8}\lambda^5,\\
\lim_{n\to \IF}b_n^2\left(A_{2n}+\frac{1}{2}\alpha\lambda^{-2}+\frac{1}{4}\lambda\right)
%&=&\lim_{n\to %\IF}b_n^2\left[\frac{b_n^2}{2}\left(\frac{1}{\lambda}-\frac{1}{\lambda_n}-\frac{\lambda_n}{2b_n^2}-\frac{3\lambda_n^3}{8b_n^4}+O(b_n^{-6})\right)
%+\frac{1}{2}\alpha\lambda^{-2}+\frac{1}{4}\lambda\right]\nonumber\\
&=&-\frac{1}{2}\beta\lambda^{-2}-\frac{1}{2}\alpha^2\lambda^{-3}+\frac{1}{4}\alpha-\frac{3}{16}\lambda^3\\
\lim_{n\to \IF}b_n^2(A_{3n}-\lambda)&=&\lim_{n\to\IF}b_n^2\left[\lambda_n+\frac{\lambda^3}{2b_n^2}+O(b_n^{-4})-\lambda\right]=-\alpha+\frac{1}{2}\lambda^3.
\end{eqnarray*}
\aH{Hence, using further \eqref{eq3.10}, \eqref{eq3.11} and \eqref{eq3.9} we obtain}
\begin{eqnarray}\label{eq3.18}
&&b_n^2\left[b_n^2\int_y^{\IF}\left(\lambda+\frac{x-z}{2\lambda}-\frac{u_n(x)-\rho_nu_n(z)}{\sqrt{1-\rho_n^2}}\right)
\varphi\left(\lambda+\frac{x-z}{2\lambda}\right)e^{-z} \d z -\kappa_1(x,y,\lambda,\alpha)\right]\nonumber\\
&=&b_n^2\left[A_{1n}+A_{2n}x-\left(\alpha-\frac{1}{2}\lambda^3-\frac{1}{2}\alpha  \lambda^{-2}x-\frac{1}{4}\lambda x\right)\right]
\EH{I_0} %\int_y^{\IF}\varphi\left(\lambda+\frac{x-z}{2\lambda}\right)e^{-z} \d z
\nonumber\\
&&-b_n^2\left[A_{2n}+A_{3n}-\left(\frac{3}{4}\lambda-\frac{1}{2}\alpha \lambda^{-2}\right)\right]
\EH{I_1} %\int_y^{\IF}\varphi\left(\lambda+\frac{x-z}{2\lambda}\right)e^{-z} z \d z
\nonumber\\
&\to&\left(\frac{1}{2}\lambda^6+2\alpha \lambda x -\lambda^4 x-2\alpha^2\right)
e^{-x}\left(1-\Phi\left(\lambda+\frac{y-x}{2\lambda}\right)\right)\nonumber\\
&&+\left(2\beta+2\alpha^2\lambda^{-1}+3\alpha\lambda^2-\frac{5}{4}\lambda^5\right)
e^{-x}\varphi\left(\lambda+\frac{y-x}{2\lambda}\right)
\end{eqnarray}
as $n \to \IF$. Consequently, using \eqref{I2}, \eqref{I3}, \eqref{eq3.9} and combining with \aH{the limits of $A_{in},i=1,2,3$} we have
\begin{eqnarray}\label{eq3.19}
&&\frac{1}{2}b_n^4\int_y^{\IF}v_n^2(x,z,\lambda)
\left(\lambda+\frac{x-z}{2\lambda}\right)
\varphi\left(\lambda+\frac{x-z}{2\lambda}\right)e^{-z} \d z\nonumber\\
&=&(A_{1n}+A_{2n}x)^2\left(\frac{\lambda}{2}+\frac{x}{4\lambda}\right)\EH{I_0}
%\int_y^{\IF}\varphi\left(\lambda+\frac{x-z}{2\lambda}\right)e^{-z} \d z
\nonumber\\
&&-\left[(A_{1n}+A_{2n}x)^2\frac{1}{4\lambda}+(A_{1n}+A_{2n}x)(A_{2n}+A_{3n})\left(\lambda+\frac{x}{2\lambda}\right)\right]
\EH{I_1}
%\int_y^{\IF}\varphi\left(\lambda+\frac{x-z}{2\lambda}\right)e^{-z} z \d z
\nonumber\\
&&+\left[(A_{1n}+A_{2n}x)(A_{2n}+A_{3n})\frac{1}{2\lambda}+(A_{2n}+A_{3n})^2\left(\frac{\lambda}{2}+\frac{x}{4\lambda}\right)\right]
\EH{I_2}
%\int_y^{\IF} \varphi\left(\lambda+\frac{x-z}{2\lambda}\right)e^{-z} z^2 \d z
\nonumber\\
&&-(A_{2n}+A_{3n})^2\frac{1}{4\lambda}
\zw{I_3}
%\int_y^{\IF}\varphi\left(\lambda+\frac{x-z}{2\lambda}\right)e^{-z} z^3 \d z
\nonumber\\
&\to&
%\left(\alpha-\frac{\lambda^3}{2}-\frac{\alpha
%x}{2\lambda^{2}}-\frac{\lambda x}{4}\right)^2
%\left(\frac{\lambda}{2}+\frac{x}{4\lambda}\right)\int_y^{\IF}
%\varphi\left(\lambda+\frac{x-z}{2\lambda}\right)e^{-z} \d z\nonumber\\
%&&-\left[\peng{\frac{1}{4\lambda}\left(\alpha-\frac{\lambda^3}{2}-\frac{\alpha
%x}{2\lambda^{2}}-\frac{\lambda x}{4}\right)^{2}
%+\left(\alpha-\frac{\lambda^3}{2}-\frac{\alpha
%x}{2\lambda^{2}}-\frac{\lambda x}{4}\right)}
%\left(\frac{3\lambda}{4}-\frac{\alpha}{2\lambda^{2}}\right)\left(\lambda+\frac{x}{2\lambda}\right)\right]
%\int_y^{\IF}
%\varphi\left(\lambda+\frac{x-z}{2\lambda}\right)e^{-z} z \d z\nonumber\\
%&&+\left[\left(\alpha-\frac{\lambda^3}{2}-\frac{\alpha
%x}{2\lambda^{2}}-\peng{\frac{\lambda x}{4}}\right)
%\left(\frac{3\lambda}{4}-\frac{\alpha}{2\lambda^{2}}\right)\frac{1}{2\lambda}
%+\left(\frac{3\lambda}{4}-\frac{\alpha}{2\lambda^{2}}\right)^2\left(\frac{\lambda}{2}+\frac{x}{4\lambda}\right)\right]
%\int_y^{\IF}
%\varphi\left(\lambda+\frac{x-z}{2\lambda}\right)e^{-z} z^2 \d z\nonumber\\
%&&-\left(\frac{3\lambda}{4}-\frac{\alpha}{2\lambda^{2}}\right)^2\frac{1}{4\lambda}
%\int_y^{\IF}
%\varphi\left(\lambda+\frac{x-z}{2\lambda}\right)e^{-z} z^3 \d z\nonumber\\
%&=&
\left(2\lambda^8+\frac{15}{2}\lambda^6-8\alpha\lambda^3+2\alpha^2-3\lambda^4x-4\lambda^6x+2\alpha
\lambda x+2\lambda^4x^2\right)
e^{-x}\overline \Phi\left(\lambda+\frac{y-x}{2\lambda}\right)\nonumber\\
&&+\left(-\frac{3}{8}\lambda^3xy-\alpha \lambda^2x+\frac{3}{4}\alpha
y^2 -\frac{\alpha^2
y^2}{4\lambda^3}-\frac{\alpha^2x^2}{4\lambda^3}+\frac{\alpha^2
xy}{2\lambda^3}
-\alpha \lambda^2y-\frac{1}{4}\alpha x^2-\frac{7}{4}\lambda^7+\frac{7}{2}\lambda^5x-\frac{1}{16}\lambda^3x^2-\frac{9}{2}\lambda^5\right.\nonumber\\
&&\left.-\alpha\lambda^4+\alpha^2\lambda+\frac{3}{2}\lambda^5y-\frac{9}{16}\lambda^3y^2+6\alpha\lambda^2-\frac{2\alpha^2}{\lambda}-\frac{1}{2}\alpha
xy\right) e^{-x}\varphi\left(\lambda+\frac{y-x}{2\lambda}\right), n\to\IF
\end{eqnarray}
where $v_n(x,z,\lambda)$ is as in the previous lemma.  Using Taylor expansion with Lagrange remainder term, we have 
\begin{eqnarray}\label{eq3.20}
\Phi\left(\frac{u_n(x)-\rho_n u_n(z)}{\sqrt{1-\rho_n^2}}\right)
&=&\Phi\left(\lambda+\frac{x-z}{2\lambda}\right)+
\varphi\left(\lambda+\frac{x-z}{2\lambda}\right)
v_n(x,y,\lambda) \Bigl[1 - \frac{1}{2} \left(\lambda+\frac{x-z}{2\lambda}\right)
v_n(x,y,\lambda)\Bigr]\notag \\ %\left(\frac{u_n(x)-\rho_n u_n(z)}{\sqrt{1-\rho_n^2}}-\lambda-\frac{x-z}{2\lambda}\right)\nonumber\\
%&&-\frac{1}{2}\varphi\left(\lambda+\frac{x-z}{2\lambda}\right)\left(\lambda+\frac{x-z}{2\lambda}\right)
%v_n^2(x,y,\lambda)\notag \\ %\left(\frac{u_n(x)-\rho_n u_n(z)}{\sqrt{1-\rho_n^2}}-\lambda-\frac{x-z}{2\lambda}\right)^2\nonumber\\
&&+\frac{1}{6}\varphi(\EH{\xi_{n}(x,z)})(\EH{\xi_{n}^2(x,z)} -1) v_n^3(x,y,\lambda)
%\left(\frac{u_n(x)-\rho_n u_n(z)}{\sqrt{1-\rho_n^2}}-\lambda-\frac{x-z}{2\lambda}\right)^3
\end{eqnarray}
for some \peng{$\EH{\xi_{n}(x,z)}$} between $\frac{u_n(x)-\rho_n
u_n(z)}{\sqrt{1-\rho_n^2}}$ and $\lambda+\frac{x-z}{2\lambda}$.
\aH{Since further}
\begin{equation}\label{eq3.21}
b_n^4\int_y^{\IF}\left(\frac{u_n(x)-\rho_n
u_n(z)}{\sqrt{1-\rho_n^2}}-\lambda-\frac{x-z}{2\lambda}\right)^3
(\EH{\xi_{n}^2(x,z)}-1)\varphi(\EH{\xi_{n}(x,z)}) e^{-z} \d z=O(b_n^{-2})
\end{equation}
\peng{the desired result follows by \eqref{eq3.18}-\eqref{eq3.21}}.
\QED

\prooftheo{th2.1} Define \[h_n(x,y,\lambda)=n\ln
F_{\rho_n}(u_n(x),u_n(y))+\Phi\left(\lambda+\frac{x-y}{2\lambda}\right)e^{-y}
+\Phi\left(\lambda+\frac{y-x}{2\lambda}\right)e^{-x}.\] \peng{\aH{In view of} \eqref{eq1.3}, \eqref{eq3.5} and Lemma \ref{le3.2}} we have
\begin{eqnarray*}
b_n^2h_n(x,y,\lambda)&=&b_n^2\left[-n(1-F_{\rho_n}(u_n(x),u_n(y)))-\frac{n}{2}(1-F_{\rho_n}(u_n(x),u_n(y)))^2(1+o(1))\right.\\
&& \left.+\Phi\left(\lambda+\frac{x-y}{2\lambda}\right)e^{-y}
+\Phi\left(\lambda+\frac{y-x}{2\lambda}\right)e^{-x}\right]\\
%&=&b_n^2[-n(1-\Phi(u_n(x)))+e^{-x}]
%+b_n^2\int_y^{\IF}\left(\Phi\left(\lambda+\frac{x-z}{2\lambda}\right)-\Phi\left(\frac{u_n(x)-\rho_nu_n(z)}{\sqrt{1-\rho_n^2}}\right)\right)e^{-z}\d %z\\
%&&-\int_y^{\IF}\Phi\left(\frac{u_n(x)-\rho_nu_n(z)}{\sqrt{1-\rho_n^2}}\right)e^{-z}\left(1-\frac{z^2}{2}\right)\d z
%+\peng{O(b_n^{-2})}-\frac{b_n^2n}{2}(1-F_{\rho_n}(u_n(x),u_n(y)))^2(1+o(1))\\
&\to&\left(\frac{1}{2}x^2+x\right)e^{-x}+\kappa_1(x,y,\lambda,\alpha)
- %\EH{[I_0 - I_2/2]}
\zw{\int_y^{\IF}\Phi\left(\lambda+\frac{x-z}{2\lambda}\right)e^{-z}\left(1-\frac{z^2}{2}\right)\d z}
\end{eqnarray*}
as $n \to \IF$.
\zw{
By partial integration we have
\begin{eqnarray*}
&&\int_y^{\IF}\Phi\left(\lambda+\frac{x-z}{2\lambda}\right)e^{-z}\left(1-\frac{z^2}{2}\right)\d z\\
&=&-\left(\frac{1}{2}y^2+y\right)e^{-y}\Phi\left(\lambda+\frac{x-y}{2\lambda}\right)
+\left(2\lambda^4-2\lambda^2x +x+\frac{1}{2}x^2\right)e^{-x}\peng{\overline \Phi\left(\lambda+\frac{y-x}{2\lambda}\right)}\\
&&+(-2\lambda^3+\lambda x+\lambda
y+2\lambda)e^{-x}\varphi\left(\lambda+\frac{y-x}{2\lambda}\right).
\end{eqnarray*}
}
Further  as $n\to \IF$
\begin{equation}\label{eq3.22}
h_n(x,y,\lambda)\to 0 \quad
\mbox{   and   }\quad
\left|\sum_{i=2}^{\IF}\frac{h_n^{i-2}(x,y,\lambda)}{i!}\right|<\exp(h_n(x,y,\lambda))\to 1.
\end{equation}
Hence\peng{,}
\begin{eqnarray*}
b_n^2\Big(F_{\rho_n}^n(u_n(x),u_n(y))-H_{\lambda}(x,y)\Big)&=&b_n^2\Big(\exp(h_n(x,y,\lambda))-1\Big)H_{\lambda}(x,y)\\
&=&b_n^2h_n(x,y,\lambda)\left(1+h_n(x,y,\lambda)\sum_{i=2}^{\IF}\frac{h_n^{i-2}(x,y,\lambda)}{i!}\right)H_{\lambda}(x,y)\\
&\to&\kapx H_{\lambda}(x,y)
\end{eqnarray*}
as $n \to \IF$,
where
\begin{eqnarray*}
\kapx &=&2^{-1}(x^2+2x)e^{-x}\Phi\left(\lambda+\frac{y-x}{2\lambda}\right)
+2^{-1}(y^2+2y)e^{-y}\Phi\left(\lambda+\frac{x-y}{2\lambda}\right)\\
&&+(2\alpha-\lambda^3-\lambda x-\lambda y-2\lambda)e^{-x}\varphi\left(\lambda+\frac{y-x}{2\lambda}\right).
\end{eqnarray*}
The proof is complete. \QED

\prooftheo{th2.2} \peng{By arguments similar to that of} Lemma
\ref{le3.2}, we have
\begin{eqnarray}\label{addpeng1}
&&b_n^2\int_y^{\IF}\left(\Phi\left(\lambda+\frac{x-z}{2\lambda}\right)-\Phi\left(\frac{u_n(x)-\rho_nu_n(z)}{\sqrt{1-\rho_n^2}}\right)\right)
e^{-z}\left(1-\frac{z^2}{2}\right) \d z \nonumber\\
&=&(A_{1n}+A_{2n}x)\zw{\left(I_0-\frac{1}{2}I_2\right)}
%\int_y^{\IF}\varphi\left(\lambda+\frac{x-z}{2\lambda}\right)e^{-z}\left(1-\frac{z^2}{2}\right)
%\d z\nonumber\\
%&&
-(A_{2n}+A_{3n})\zw{\left(I_1-\frac{1}{2}I_3\right)}
%\int_y^{\IF}\varphi\left(\lambda+\frac{x-z}{2\lambda}\right)e^{-z}z\left(1-\frac{z^2}{2}\right) \d z
\peng{+O(b_{n}^{-2})}\nonumber\\
&\to&\left(\alpha-\frac{\lambda^3}{2}-\frac{\alpha
x}{2\lambda^{2}}-\frac{\lambda x}{4}\right)\zw{\left(I_0-\frac{1}{2}I_2\right)}
%\int_y^{\IF}\left(1-\frac{z^2}{2}\right)e^{-z}\peng{\varphi\left(\lambda+\frac{x-z}{2\lambda}\right)} \d z \nonumber\\
-\left(\frac{3}{4}\lambda-\frac{\alpha}{2 \lambda^2}\right)\zw{\left(I_1-\frac{1}{2}I_3\right)}
%\int_y^{\IF}z\left(1-\frac{z^2}{2}\right) \peng{e^{-z}\varphi\left(\lambda+\frac{x-z}{2\lambda}\right)}\d z
\nonumber\\
&=&\peng{(2\lambda^4-4\alpha \lambda x-2\lambda^2 x+8\lambda^6
x-5\lambda^4 x^2+10\lambda^4 x+\lambda^2
x^3+8\alpha\lambda^3-4\lambda^8-16\lambda^6)
e^{-x}\overline \Phi\left(\lambda+\frac{y-x}{2\lambda}\right)}\nonumber\\
&&+\peng{\left(2\alpha+4\lambda^7+12\lambda^5-3\lambda^3-6\lambda^5x+2\lambda^3x^2-\alpha
y^2+2\lambda^3xy-2\lambda^5y+\frac{3}{2}\lambda^3y^2-8\alpha
\lambda^2\right)
e^{-x}\varphi\left(\lambda+\frac{y-x}{2\lambda}\right)}\nonumber\\
&=&\tau_2(\alpha, \lambda,x,y)
\end{eqnarray}
as $n \to \IF$.
By partial integration we get
\begin{eqnarray}\label{addpeng2}
&&\int_y^{\IF}\Phi\left(\lambda+\frac{x-z}{2\lambda}\right)e^{-z}\left(\frac{z^4}{8}-\frac{z^2}{2}-\peng{2}\right)\d z\nonumber\\
&=&8^{-1}(y^4+4y^3+8y^2+16y)e^{-y}\Phi\left(\lambda+\frac{x-y}{2\lambda}\right)
-\frac{1}{16\lambda}\int_y^{\IF}\varphi\left(\lambda+\frac{x-z}{2\lambda}\right)e^{-z}(z^4+4z^3+8z^2+16z) \d z\nonumber\\
&=&8^{-1}(y^4+4y^3+8y^2+16y)e^{-y}\Phi\left(\lambda+\frac{x-y}{2\lambda}\right)\nonumber\\
&&+\Big(4
\lambda^6x-3\lambda^4x^2+\lambda^2x^3-2\lambda^8-8\lambda^6-\frac{1}{8}x^4-2\lambda^2x-\frac{1}{2}x^3+\peng{2}\lambda^4+6\lambda^4x-x^2-2x\Big)
e^{-x}\overline \Phi\left(\lambda+\frac{y-x}{2\lambda}\right)\nonumber\\
&&+\Big(2\lambda^7-\frac{1}{4}\lambda
x^3-\peng{\lambda^5y+\frac{1}{2}\lambda^3y^2+\lambda^3xy-3\lambda^5x+\frac{3}{2}\lambda^3x^2}-\frac{1}{4}\lambda
y^3
-\frac{1}{4}\lambda y^2x-\frac{1}{4}\lambda yx^2-\lambda^3y-\lambda y^2-\lambda xy-\lambda^3x\nonumber\\
&&-\lambda x^2-4\lambda^3+6\lambda^5-2\lambda x-2\lambda
y-4\lambda\Big)
e^{-x}\varphi\left(\lambda+\frac{y-x}{2\lambda}\right)\nonumber\\
&=&\tau_3(\lambda, x,y).
\end{eqnarray}
Hence, using \eqref{eq1.3}, \eqref{eq3.4} and Lemma \ref{le3.3} we have
\begin{eqnarray}\label{addpeng3}
&&b_n^2\Big[b_n^2h_n(\lambda,x,y)-\kapx \Big]\nonumber\\
&=&b_n^2\Biggl( b_n^2\Biggl( -n(1-F_{\rho_n}(u_n(x),u_n(y)))-\frac{n}{2}(1-F_{\rho_n}(u_n(x),u_n(y)))^2(1+o(1))
+\Phi\left(\lambda+\frac{x-y}{2\lambda}\right)e^{-y} \nonumber\\
&&+ \Phi\left(\lambda+\frac{y-x}{2\lambda}\right)e^{-x}\Biggr)-\kapx \Biggr)\nonumber\\
&=&b_n^2\Big[b_n^2\left[-\pHn(x)+e^{-x}\right]-\left(\frac{x^2}{2}+x\right)e^{-x}\Big]\nonumber\\
&&+b_n^2\left[b_n^2\int_y^{\IF}\left(\Phi\left(\lambda+\frac{x-z}{2\lambda}\right)
-\Phi\left(\frac{u_n(x)-\rho u_n(z)}{\sqrt{1-\rho^2}}\right)\right)e^{-z} \d z -\kappa_1(x,y,\lambda,\alpha)\right]\nonumber\\
&&+b_n^2\int_y^{\IF}\left(\Phi\left(\lambda+\frac{x-z}{2\lambda}\right)-\Phi\left(\frac{u_n(x)-\rho_nu_n(z)}{\sqrt{1-\rho_n^2}}\right)\right)
e^{-z}\left(1-\frac{z^2}{2}\right) \d z \nonumber\\
&&-\int_y^{\IF}\Phi\left(\lambda+\frac{x-z}{2\lambda}\right)e^{-z}\left(\frac{z^4}{8}-\frac{z^2}{2}-\peng{2}\right)\d
z \wz{+O(b_n^{-2})}
-\frac{1}{2}b_n^4n(1-F_{\rho_n}(u_n(x),u_n(y)))^2(1+o(1))\nonumber\\
&\to&-8^{-1}(x^4+4x^3+8x^2+16x)e^{-x}+\tau_1(\alpha,\beta,\lambda, x,y)+\tau_2(\alpha,\lambda, x,y)-\tau_3(\lambda,x,y)\nonumber\\
&=&\tapx, \quad n\to \IF, 
\end{eqnarray}
where \EH{$\tau_i,i\le 3$} are given by Lemma \ref{le3.3}, \eqref{addpeng1} and \eqref{addpeng2}, respectively. Hence,
\eqref{eq3.22} entails
\begin{eqnarray*}
&&b_n^2\Big[b_n^2\Big(F_{\rho_n}^n(u_n(x),u_n(y))-H_{\lambda}(x,y)\Big)-\kapx H_{\lambda}(x,y)\Big]\\
&=&b_n^2\Big[b_n^2\Big(\exp(h_n(x,y,\lambda))-1\Big)-\kapx \Big]H_{\lambda}(x,y)\\
&=&\left[b_n^2\Big[b_n^2h_n(x,y,\lambda)-\kapx \Big]
+b_n^4h_n^2(x,y,\lambda)\left(\frac{1}{2}+h_n(x,y,\lambda)\sum_{i=3}^{\IF}\frac{h_n^{i-3}(x,y,\lambda)}{i!}\right)\right]H_{\lambda}(x,y)\\
&\to&\left(\tapx +\frac{1}{2}\kappa^2(\alpha, \lambda,x,y)\right)H_{\lambda}(x,y)
\end{eqnarray*}
as $n\to \IF$ establishing the proof. \QED

\prooftheo{th2.3} \pzx{$(i)$ For the case of $\rho_{n}\in [-1, 0]$,
we first consider that the bivariate Gaussian are either complete
independent ($\rho_{n}=0$) or complete negative dependent
($\rho_{n}=-1$). Both imply $\lambda=\infty$. Let}
$\hat{h}_n(x,y)=n\ln F_0(u_n(x),u_n(y))+e^{-x}+e^{-y}$ and
$\tilde{h}_n(x,y)=n\ln F_{-1}(u_n(x),u_n(y))+e^{-x}+e^{-y}$.
In view of Lemma 2.1 in \cite{MR606804}
\begin{eqnarray*}
b_n^2\hat{h}_n(x,y)&=&b_n^2[-\pHn(x)+e^{-x}]+b_n^2[-\pHn(y)+e^{-y}]\\
&&+\frac{b_n^2}{n}\pHn(x) \pHn(y)-\frac{1}{2}b_n^2n(1-F_0(u_n(x),u_n(y)))^2(1+o(1))\\
&\to&s(x)+s(y)
\end{eqnarray*}
and
\begin{eqnarray*}
b_n^2\tilde{h}_n(x,y)&=&b_n^2[-\pHn(x)+e^{-x}]+b_n^2[-\pHn(y)+e^{-y}]\\
&&\zw{+b_n^2n\pk{u_n(x)<X<-u_n(y)}}-\frac{1}{2}b_n^2n(1-F_{-1}(u_n(x),u_n(y)))^2(1+o(1))\\
&\to&s(x)+s(y)
\end{eqnarray*}
 \peng{as $n\to\infty$, \zw{where $X$ is a standard normal variable}. \aH{By}  Lemma 2.1 in \cite{MR606804}
once} again we have
$$\lim_{n\to \IF}b_n^2\Big[b_n^2\hat{h}_n(x,y)-(s(x)+s(y))\Big]=t(x)+t(y)$$
and
$$\lim_{n\to \IF}b_n^2\Big[b_n^2\tilde{h}_n(x,y)-(s(x)+s(y))\Big]=t(x)+t(y).$$
\pzx{\aH{Consequently},
\begin{eqnarray}\label{addpeng4}
\lim_{n\to
\IF}b_n^2\Big[b_{n}^{2}\DFxyD-(s(x)+s(y))H_{\IF}(x,y)\Big]
=\left(t(x)+t(y)+\frac{1}{2}(s(x)+s(y))^2\right)H_{\IF}(x,y)
\end{eqnarray}
holds with $\rho_{n}=-1$ and $\rho_{n}=0$ (for all $n$ large), respectively.
Consequently, by using Slepian's Lemma and \eqref{addpeng4}, the claimed
result \eqref{eq2.5} holds for $\rho_{n}\in [-1,0]$.}

\zw{$(ii)$} For the complete positive \Ea{dependence}  case, without loss of generality, assume that $x<y$. \EH{Hence}
$$F_1(u_n(x),u_n(y))=\Phi(u_n(x)), \quad H_0(x,y)=\Lambda(x)$$
\eqref{eq1.3} follows and thus the proof is complete. \QED

\proofkorr{cor2.1} \EH{$(i)$} \peng{Obviously, $n(1-\Phi(b_{n}))=1$ implies
\[b_{n}\sim (2\ln  n)^{1/2},\quad e^{-\frac{b_{n}^{2}}{2}}\sim
\sqrt{2\pi}\frac{b_{n}}{n}\] for large $n$. For the case of
$\lambda=\infty$, according} to Berman's inequality \aH{(see e.g.,
\cite{Pit96})}, with some \EH{positive} constant \EH{which may change from line to line} $\mathbb{C}$ \EH{and all $n$ large} we have
\begin{eqnarray*}
\lefteqn{b_{n}^{2}\Big|F_{\rho_n}^n(u_n(x),u_n(y))-F_0^n(u_n(x),u_n(y))\Big|}\\
 &\le& \mathbb{C} n(\ln  n)
\rho_n\exp\left(-\frac{u^2_n(x)+u^2_n(y)}{2(1+\rho_n)}\right)\\
&\le& \mathbb{C} n(\ln  n)
\left(\exp\left(-\frac{u^2_n(x)}{(1+\rho_n)}\right)+
\exp\left(-\frac{u^2_n(y)}{(1+\rho_n)}\right)\right)\\
&\le &
\mathbb{C}\Big(\exp(2|x|)+\exp(2|y|)\Big)n^{-\frac{1-\rho_{n}}{1+\rho_{n}}}\left(\ln
n\right)^{1+\frac{1}{1+\rho_{n}}}\\
&\le & \mathbb{C}\Big(\exp(2|x|)+\exp(2|y|)\Big)
\exp\Big(-\frac{1}{2}\left((1-\rho_{n})\ln  n
-(2+\rho_{n})\ln \ln  n\right)\Big)\\
&\to& \EH{\mathbb{C} e^{-\gamma} \Big(\exp(2|x|)+\exp(2|y|)\Big)}, \quad n\to\infty
\end{eqnarray*}
since  by the assumption $\lim_{n\to \IF}\zw{\frac{1}{2}} ((1-\rho_{n})\ln  n
-(2+\rho_{n})\ln \ln  n)=\gamma$. \peng{By Theorem
\ref{th2.3}} we have
\begin{eqnarray*}
\limsup_{n\to\IF}b_n^2 \DFxyD %\Big|F_{\rho_n}^n(u_n(x),u_n(y))-H_{\IF}(x,y)\Big|\\
&\le&\lim_{n\to\IF}b_n^2\DFxyB %\Big|F_{0}^n(u_n(x),u_n(y))-H_{\IF}(x,y)\Big|
+\lim_{n\to\IF}b_n^2\Big|F_{\rho_n}^n(u_n(x),u_n(y))-F_{0}^n(u_n(x),u_n(y))\Big|\\
&\le &(|s(x)|+|s(y)|)H_{\IF}(x,y)+ \EH{\mathbb{C} e^{-\gamma} \Big(\exp(2|x|)+\exp(2|y|)\Big)}.
\end{eqnarray*}
\EH{$(ii)$} {The condition $\lim_{n\to\infty}(1-\rho_{n})(\ln
n)^{3}=\PZX{\tau^{2}\in[0,\infty)}$ implies $\lambda=0$ and
$\lim_{n\to\infty}\rho_{n}=1$. By} using \aH{Berman's inequality
\wz{in \cite{MR1902188}}} we have
\begin{eqnarray*}
b_{n}^{2}\Big|F_{\rho_n}^n(u_n(x),u_n(y))-F_1^n(u_n(x),u_n(y))\Big|
&\le& \mathbb{C} n(\ln  n)
\wz{\left(\frac{\pi}{2}-\arcsin(\rho_n)\right)}\exp\left(-\frac{u^2_n(x)+u^2_n(y)}{\wz{4}}\right)\\
&\le&\mathbb{C}\Big(\exp(2|x|)+\exp(2|y|)\Big)\left(\frac{\pi}{2}-\arcsin(\rho_n)\right)(\ln  n)^{\frac{3}{2}}\\
&\le& \mathbb{C}\Big(\exp(2|x|)+\exp(2|y|)\Big)\wz{\left(1-\rho_n \right)^{\frac{1}{2}}\left(\ln  n\right)^{\frac{3}{2}}}\\
 &\to& \PZX{\tau} \mathbb{C}\Big(\exp(2|x|)+\exp(2|y|)\PZX{\Big)}, \quad \EH{n\to \IF}
\end{eqnarray*}
since  $\lim_{n\to \IF}\wz{\left(1-\rho_n\right)(\ln  n)^3
\EH{=\tau^2}}$ \PZX{ which also implies
$\lim_{n\to\infty}\frac{\frac{\pi}{2}-\arcsin(\rho_n)}{(1-\rho_{n})^{\wz{1/2}}}=\wz{\sqrt{2}}$}.
\wz{Hence  Theorem \ref{th2.3} \Ha{yields}}
\begin{eqnarray*}
\limsup_{n\to\IF}b_n^2 \DFxyE % \Big|F_{\rho_n}^n(u_n(x),u_n(y))-H_{0}(x,y)\Big|\\
&\le&\EH{\lim_{n\to\IF}b_n^2}\DFxyA %\lim_{n\to\IF}b_n^2\Big|F_{1}^n(u_n(x),u_n(y))-H_{0}(x,y)\Big|
+\lim_{n\to\IF}b_n^2\Big|F_{\rho_n}^n(u_n(x),u_n(y))-F_{1}^n(u_n(x),u_n(y))\Big|\\
&=&|s(\min(x,y))|H_{0}(x,y)\PZX{+\tau
\mathbb{C}\Big(\exp(2|x|)+\exp(2|y|)\Big)}
\end{eqnarray*}
establishing the claim. \QED

%\section{Numerical analysis}
\vspace{0.5in}

\noindent\textbf{Acknowledgments}. \EH{The authors are in debt to the referees and the  Associate  Editor for several suggestions that lead to various improvements. E.\ Hashorva and Z. Weng acknowledge support from}
the Swiss National Science Foundation grants  200021-140633/1, 200021-134785 and the project RARE -318984 (an FP7  Marie Curie
IRSES Fellowship);
 Z. Peng  has been supported by the National Natural Science Foundation of China under grant
11171275 and the Natural Science Foundation Project of CQ under cstc2012jjA00029.

\bibliographystyle{plain}
\bibliography{RateC}
\end{document}